\documentclass{amsart}

\theoremstyle{definition}

\theoremstyle{remark}

\numberwithin{equation}{section}

\begin{document}

\title{Realizations of Weyl groups on ellipsoids}
\author{Anatoli Loutsiouk}
\address{Dept. Math., King Mongkut University of Technology Thonburi,Thailand}
\email{loutsiouk.ana@kmutt.ac.th}
\begin{abstract}

For any finite-dimensional complex semisimple Lie algebra two ellipsoids (primary and secondary) 
    are considered. The equations of these ellipsoids are Diopantine equations and the Weyl group acts on the sets of all their Diophantine solutions. This provides two realizatons (primary and secondary) of the Weyl group. The primary realization suggests an order on the Weyl group, which is useful for an exploration of the Bruhat order. The secondary realization is useful for finding for any element of the Weyl group all of its reduced expressions, which implies another realization of the Bruhat order. 

\end{abstract}

\maketitle

\section{Introduction}
In this paper, two geometric objects are explored for any complex semisimple Lie algebra $\mathfrak{G}$, which provide two realizations for the Weyl group attached to this algebra. These objects are ellipsoids in the real linear space $\mathbb{R}^n$.  

Given a complex semisimple Lie algebra $\mathfrak{G}$ of rank $n$ and a Cartan subalgebra $\mathfrak{H}$, the pair ($\mathfrak{G}$, $\mathfrak{H}$) determines the system of roots $\Delta$, a subsystem $\Delta^+$ of all positive roots, and the subsystem $\Pi$ of all simple roots. In the real vector space $\mathbb{R}^n$, define an inner product $<\  ,\  >$ so that $2<e_i,e_j>/<e_i,e_i>=a_{ij},$ 
where $a_{ij}$ is the element of the Cartan matrix $A$ defined by the system $\Pi$ of simple roots, and $e_i$ is the element of the standard basis of $\mathbb{R}^n$, which we identify with the simple root $\alpha_i$ from $\Pi$. We also partially order the set $\mathbb{R}^n$ by setting $x \leq y$ if $x_i \leq y_i$ for all $i$.
Let $\delta$ be the half-sum of all positive roots in this realization of $\Delta^+$. The element $\delta \in \mathbb{R}^n_+$ satisfies the equation 
$A\delta=\widetilde{1}=(1,...,1).$
\section{Primary and Secondary ellipsoids}
The principal object of study in this paper is the subset of $\mathbb{R}^n$ defined by the equation
\begin{equation}<x,x-2\delta>=0.\label{elprim1} \end{equation}
We shall call this subset the {\em primary ellipsoid} and denote it by $PE(\mathfrak{G})$.
It is easy to write down Equation \ref{elprim1} in coordinate form through the Dynkin diagram for the semisimple Lie algebra $\mathfrak{G}$. The Dynkin diagram has $n$ vertices $a_i$. Each vertice  has a weight denoted by $k_i$, which is an integer 1,2,or 3.  Some of the vertices are joined by links. For any link joining verices $a_i$ and $a_j$, set $l_{ij}=max{(k_i,k_j)}$, otherwise set $l_{ij}=0$. 
\newtheorem{Theorem}{Theorem}[section]
\begin{Theorem} \label{pro:1thm1}
Equation \ref{elprim1} in coordinate form is as follows:
\begin{equation}\sum_{i=1}^{n}k_i(x_{i}^{2}-{x_i})-\sum_{1\leq i<j\leq n}{l_{ij}}{x_i}{x_j}=0,\label{elprim2} \end{equation}
where the first sum is taken over all the vertices, and the second sum is taken over all the links in the Dynkin diagram for the complex semisimple Lie algebra $\mathfrak{G}$.  
\end{Theorem}

For any root $\alpha$, let  
\begin{equation}m_{\alpha}=\frac{2<{\alpha,\delta}>}{<{\alpha,\alpha}>}.\label{elprim3} \end{equation}
\newtheorem{Proposition}{Proposition}[section]
\begin{Proposition} \label{pro:1thm2}$m_\alpha$ is an integer, which is positive if and only if the root $\alpha$ is positive, and it is equal to 1 if and only if the positive root  $\alpha$ is simple.
\end{Proposition}
We shall call the number $m_{\alpha}$ the {\em grade} of the positive root $\alpha$.

\begin{Proposition} \label{pro:1thm7} For any positive root $\alpha$=$(a_1,a_2,...,a_n)$, the element
$m_{\alpha}{\alpha}$ 
belongs to the primary ellipsoid $E(\mathfrak{G})$ defined by Equation \ref{elprim1} or \ref{elprim2}.
\end{Proposition}

Another ellipsoid related to the semisimple Lie algebra $\mathfrak{G}$ is denoted by $SE(\mathfrak{G})$, and we shall call it the {\em secondary ellipsoid} for  $\mathfrak{G}$. For any $x\in PE(\mathfrak{G})$, with $x=(x_1,...,x_n)$, and for any $i\in \{1,...,n\}$, we define $h^x_i\in \mathbb{R}\backslash \{0\}$, if it exists, otherwise we set $h^x_i = 0$, so that 
$T_i(x)=^ix=(x_1,...,x_{i-1},x_i+h^x_i, x_{i+1},...,x_n) \in PE(\mathfrak{G}).$ 
Such $h^x_i$ is a unique real number. Consider the vector $h^x=(h^x_1,h^x_2,...,h^x_n)\in \mathbb{R}^n$.

\begin{Proposition} \label{pro:1thm3}
$h^x=\widetilde{1}- Ax.$
\end{Proposition}

Observe that if $x$ is an integral vector, that is a vector with all integer components, then $h^x$ is an integral vector as well.
The set $SE(\mathfrak{G})$ is the set of all such vectors $h^x\in \mathbb{R}^n$ as $x$ runs through the primary ellipsoid $PE(\mathfrak{G})$. 
\begin{Theorem} \label{pro:1thm4}
The subset $SE(\mathfrak{G})$ is an ellipsoid, which is described by the equation
\begin{equation}<A^{-1}(\widetilde{1}-h),A^{-1}(\widetilde{1}+h)>=0,\label{elprim4}\end{equation}
\end{Theorem}
\begin{Theorem} \label{pro:1thm5} The equation of the secondary ellipsoid in coordinate form is as follows:
\begin{equation} \sum_{i=1}^{n}{k_{i}b_{ii}(h_{i}^{2}-1)}+2\sum_{1\leq i<j\leq n}{k_{i}b_{ij}(h_i}{h_j}-1)=0,\label{elprim5} \end{equation} 
where $k_{i}$ is the weight of the $i$-th vertice in the Dynkin diagram, and $b_{ii},b_{ij}$ are elements of the matrix $A^{-1}$  (the inverse matrix for the Cartan matrix $A$).
\end{Theorem}
By multiplying Equation \ref{elprim5} with $detA$, we get an equation with all coefficients being nonnegative integers.
\section{Primary and secondary realizations of the Weyl group}

The mappings $T_i$ defined in the preceding section for $1\leq i \leq n$ are involutions on the primary ellipsoid. In the group of all bijections of the primary ellipsoid, consider the subgroup $WPE(\mathfrak{G})$ generated by the mappings $T_i$. 
\begin{Theorem} \label{pro;1thm8} The group $WPE(\mathfrak{G})$ is isomorphic to the Weyl group of the Lie algebra $\mathfrak{G}$. 
\end{Theorem}

Consider Equation \ref{elprim2} of primary ellipsoid and Equation \ref{elprim5} of secondary ellipsoid as Diophantine equations. We shall denote  by $PD(\mathfrak{G})$ and $SD(\mathfrak{G})$ their respective sets of Diophantine solutions. 

\begin{Theorem} \label{pro;1thm9} The subset $PD(\mathfrak{G})$ is invariant under the action of the group $WPE(\mathfrak{G})$, and it splits into orbits. The set of the orbits is one-to-one with a subset of the set of all integral solutions with all nonnegative components of the equation of secondary ellipsoid. For any such a solution $h$, the vector $x_h=A^{-1}(\widetilde{1}-h)$, if it is integral, is the unique minimal vector of the respective orbit under the partial ordering $\leq$.
\end{Theorem}

We shall parametrise the set of the orbits by the family of their minimal vectors. For any such vector $a=x_h$, denote by $PD_a(\mathfrak{G})$ the respective orbit. The number of orbits for the case $ E_8$ is 157. 
There is only one integral solution $h$ with all nonnegative components that has all components positive, and it is equal to $\widetilde{1}$. The respective vector $x_{\widetilde{1}}=0$. All other such solutions have at least one component equal to $0$. For any such $h$ consider all those values of index $i_1<1_2<...<i_k$ for which the corresponding component is $0$. In the Weyl group $W$, let $W_h$ be the subgroup generated by the elementary reflections $e_{i_1}, e_{i_2},...,e_{i_k}$.
\begin{Theorem} \label{pro;1thm10}The number of elements in the respective orbit is equal to the number $card (W)/card(W_h)$. 
\end{Theorem}
There is only one orbit with the number of elements equal to the order of the Weyl group. This orbit contains the origin $0$ and all the  vectors $m_{\alpha}\alpha$ for any positive root $\alpha$, as well as the element $2\delta$ together with all the elements $2\delta-m_{\alpha}\alpha$. We shall call this orbit {\em the main orbit} and denote it by $PD_0(\mathfrak{G})$.
The corresponding subset of the secondary ellipsoid will be denoted by $SD_0(\mathfrak{G}).$  

The procedure to find the orbits $PD_a(\mathfrak{G})$ and $SD_a(\mathfrak{G})$ is as follows:

Start with $a$ belonging to $PD_a(\mathfrak{G})$. Consider $h^a=\widetilde{1}-Aa$ belonging to $SD_a(\mathfrak{G})$. For any positive component $h^{a}_i$ consider the vector $^{i}a=(a_1,...,a_{i}+h^{a}_i,...,a_n)$. This vector belongs to $PD_a(\mathfrak{G})$. Apply the procedure again to any vector $x\in PD_a(\mathfrak{G})$ thus obtained until the vector $x$ is such that $h^x$ has no more positive components. The result of this are the sets $PD_a(\mathfrak{G})$  and $SD_a(\mathfrak{G})$. 

\begin{Proposition} \label{pro:1thm11} The vectors from $PD_0(\mathfrak{G})$ do not have negative components and are sums of distinct positive roots, and the vectors of $SD_0(\mathfrak{G})$ do not have zero components.
\end{Proposition}
Denote by $MW(\mathfrak{G})$ the Weyl group of the complex semisimple Lie algebra $\mathfrak{G}$ realized in matrix form with respect to the basis of simple roots. The matrices of this matrix group have all their entries integers.
Assign to any element $w$ of the Weyl group $MW(\mathfrak{G})$ the vector 
$P(w)=\delta - w\delta.$
\begin{Theorem} \label{pro;1thm12} The mapping $P$ is a bijection from $MW(\mathfrak{G})$ to $PD_0(\mathfrak{G})$. 
\end{Theorem}

So, we shall consider $PW(\mathfrak{G})$ to be the set $PD_0(\mathfrak{G})$ with the transfered operation $*$, and call it the primary realization of the Weyl group.
In this realization of the Weyl group the identity element is the origin $0=(0,...,0)$. This realization has some interesting features. Here is one of them.

Given a positive root $\alpha$ and $b$ an arbitrary element of $PW(\mathfrak{G})$. Then, since $m_{\alpha}\alpha \in PW(\mathfrak{G}),$ 
$m_{\alpha}\alpha*b$ exists and belongs to $PW(\mathfrak{G})$. 

\begin{Theorem} \label{pro;1thm13} There exists a nonzero integer $p_{\alpha,b}$ such that $m_{\alpha}\alpha*b = p_{\alpha,b}\alpha+b$. 
\end{Theorem} 

Assign to any element $w$ of the Weyl group $MW(\mathfrak{G})$ the vector
$S(w)=Aw\delta.$

\begin{Theorem} \label{pro;1thm14} The mapping $S$ is a bijection from $MW(\mathfrak{G})$ to $SD_0(\mathfrak{G}).$ 
\end{Theorem}

This provides the secondary realization $SW(\mathfrak{G})$ of the Weyl group on the set $SD_0(\mathfrak{G})$ with the group operation transfered from matrix group $MW(\mathfrak{G})$ by the mapping $S$. The identity element is $\widetilde{1}$.
\section{Orderings of the Weyl group}

The realization of the Weyl group on the primary ellipsoid $PW(\mathfrak{G})$ provides a partial  ordering of this group that is inherited from the natural partial ordering of the linear space $\mathbb{R}^n$.    We call it the {\em primary ordering} of the group $PW(\mathfrak{G})$. 
There is another very important for different applications ordering for any Coxeter group, which is called Bruhat ordering, see\cite{bjorner,humphreys1,humphreys2}, and which we denote by $\prec$. 
 
\begin{Theorem} \label{pro;1thm15} On a Weyl group, $a\prec b$ implies  $a\leq b$.
\end{Theorem}
As a matter of fact, the primary ordering is in some cases strictly stronger than the Bruhat ordering on the Weyl group, as one can see, for example, from the case of $A_3$. In this case, the Weyl group is isomorphic to the symmetric group $S_4$, and the graph of the Bruhat ordering for this group is available in Fig. 2.4 of [1]. When compared to the primary ordering, it can be seen that there are two cases of discrepancy between the two orders for this Weyl group. In the case of the Bruhat ordering, the elements (1 4 3 2) and (4 1 2 3) of $S_4$ are not comparable, as one can see from Fig. 2.4 of [1], but their respective counterparts in the primary ellipsoid realization are the vectors (0,2,2) and (1,2,3), which are comparable in the primary order. The same applies to the elements (3 2 1 4) and (2 3 4 1) of $S_4$, which have the vectors (2,2,0) and (3,2,1) as their respective counterparts. In all other cases, the two orders agree for $A_3$. Observe that in these two cases of $a\leq b$ but not $a\prec b$ we have that $b-a=(1,0,1)$, which is not a multiple of a positive root.

To Bruhat order a Weyl group $W$ by using the primary realization take the following steps: 1)Realize W primarily; 2) Order the pimary realization primarily by inserting a link between any two directly adjacent elements $a\leq b$; 3)Delete all those links with $a\leq b$ for which $b-a$ is not a multiple of a positive root. The remaining links provide the Bruhat ordering of the Weyl group $W$.

The secondary realization of a Weyl group provides an efficient way to obtain all the reduced expressions for any ellement $w$ of the Weyl group. A reduced expression of $w$ is a shortest possible expression of it as a product of simple reflections $s_i$. Finding all reduced expressions of any element of a Weyl group boils down to finding the first element in any such expression, because if $s$ is known to be the first element of a reduced expresion for $w$, to find the second element of this reduced expression is equivalent to finding the first element of the product $sw$, and so on.
\begin{Theorem} \label{pro;1thm16} Given an element  $w$ of the Weyl group  $MW(\mathfrak{G})$, consider its image $S(w)=Aw\delta$ in $SW(\mathfrak{G})$. The vector $S(w)$ being an n-tuple of positive and negative integers $(b_1, b_2,...,b_n),$ let $i_1<i_2<..<i_m$ be the values of index $i$ for which $b_i$ is negative. Then $s_{i_1}, s_{i_2},..., s_{i_m}$ are the only simple reflections that can be the first elements of a reduced expression of $w$.
\end{Theorem}
This theorem provides an alternative way to build the Bruhat ordering of a Weyl group, because as is well known for any Coxeter group (see for example \cite{humphreys2}), knowing reduced expressions leads to Bruhat ordering through subexpressions.

%

\end{document}